\documentclass[journal]{IEEEtran}
\usepackage{graphicx}
\usepackage{amsmath}
\usepackage{amssymb}
\usepackage{CJK}
\usepackage{color}
\usepackage{dsfont}
\usepackage{subfigure}
\usepackage{booktabs}
\usepackage{microtype}
\usepackage{mathrsfs}
\usepackage{bm}
\usepackage{cite}
\usepackage{url}
\usepackage[colorlinks=true,allcolors=blue]{hyperref}

\ifCLASSINFOpdf
\else
\fi

\begin{document}

\title{Data-Based Distributionally Robust \\Stochastic Optimal Power Flow  -- Part I: \\ Methodologies}

\author{Yi Guo,~\IEEEmembership{Student Member,~IEEE,}
        Kyri Baker,~\IEEEmembership{Member,~IEEE,}
        Emiliano Dall'Anese,~\IEEEmembership{Member,~IEEE,}\\
        Zechun Hu,~\IEEEmembership{Senior Member,~IEEE,}
        Tyler H. Summers,~\IEEEmembership{Member,~IEEE}

\thanks{This material is based on work supported by the National Science Foundation under grant CNS-1566127. \emph{(Corresponding author: Tyler H. Summers)}}
\thanks{Y. Guo and T.H. Summers are with the Department
of Mechanical Engineering, The University of Texas at Dallas, Richardson,
TX, USA, email: \{yi.guo2,tyler.summers\}@utdallas.edu.}
\thanks{K. Baker is with the Department of Civil, Environmental, and Architectural Engineering, University of Colorado Boulder, Boulder, CO, USA, email:kyri.baker@colorado.edu.}
\thanks{E. Dall'Anese is with the Department of Electrical, Computer, and Energy Engineering, University of Colorado Boulder, Boulder, CO, USA, email: emiliano.dallanese@colorado.edu.}
\thanks{Z. Hu is with the Department of Electrical Engineering, Tsinghua University, Beijing, China, email: zechhu@tsinghua.edu.cn.}
}

\maketitle
\begin{abstract}
We propose a data-based method to solve a multi-stage stochastic optimal power flow (OPF) problem based on limited information about forecast error distributions. The framework explicitly combines multi-stage feedback policies with any forecasting method and historical forecast error data. The objective is to determine power scheduling policies for controllable devices in a power network to balance operational cost and conditional value-at-risk (CVaR) of device and network constraint violations. These decisions include both nominal power schedules and reserve policies, which specify planned reactions to forecast errors in order to accommodate fluctuating renewable energy sources. Instead of assuming the uncertainties across the networks follow prescribed probability distributions, we consider \emph{ambiguity sets} of distributions centered around a finite training dataset. By utilizing the Wasserstein metric to quantify differences between the empirical data-based distribution and the real unknown data-generating distribution, we formulate a multi-stage distributionally robust OPF problem to compute control policies that are robust to both forecast errors and sampling errors inherent in the dataset. Two specific data-based distributionally robust stochastic OPF problems are proposed for distribution networks and transmission systems.
\end{abstract}



%
\IEEEpeerreviewmaketitle

\section{Introduction}
\IEEEPARstart{T}{he} continued integration of renewable energy sources (RESs) in power systems is making it more complicated for system operators to balance economic efficiency and system reliability and security. As penetration levels of RESs reach substantial fractions of total supplied power, networks will face high operational risks under current operational paradigms. As it becomes more difficult to predict the net load, large forecast errors can lead to power security and reliability issues causing significant damage and costly outages. Future power networks will require more sophisticated methods for managing these risks, at both transmission and distribution levels.

The flexibility of controllable devices, including power-electronics-interfaced RESs, can be utilized to balance efficiency and risk with  optimal power flow methods \cite{opf,dommel,alsac,Baldick,low1,low2,madani2015convex,li2017linear,li2017distribued}, which aim to determine power schedules for controllable devices in a power network to optimize an objective function. However, most OPF methods in the research literature and those widely used in practice are deterministic, assuming point forecasts of exogenous power injections and ignoring forecast errors. Increasing forecast errors push the underlying distributed feedback controllers that must handle the transients caused by these errors closer to stability limits \cite{kundur}.

More recently, research focus has turned to stochastic and robust optimal power flow methods that explicitly incorporate forecast errors, in order to more systematically trade off economic efficiency and risk and to ease the burden on feedback controllers \cite{yong,capitanescu,conejo,bienstock,Vrakopoulou1,vrakopoulou2,warrington,zhang1,perninge,roald1,tyler,jabr,roald2,roald3,LiMathieu,zhang2,lubin,Baker,DallAnese1,vrakopoulou3
,li2017chance,louca2016stochastic,venzke2018convex,xie2018distributionally} Many formulations assume that uncertain forecast errors follow a prescribed probability distribution (commonly, Gaussian \cite{bienstock,lubin,roald1,venzke2018convex}) and utilize analytically tractable reformulations of probabilistic constraints. However, such assumptions are unjustifiable due to increasingly complex, nonlinear phenomena in emerging power networks, and can significantly underestimate risk.

In practice, forecast error probability distributions are never known; they are only observed indirectly through finite datasets. Sampling-based methods have been applied with a focus on quantifying the probability of constraint violation \cite{Vrakopoulou1,vrakopoulou2} and for constraining or optimizing conditional value at risk (CVaR) \cite{zhang1,tyler,DallAnese1}. 
The prediction-realization approach \cite{mohagheghi2017framework,mohagheghi2018real} solves an online stochastic optimal power flow problem by a reconciliation algorithm, which ensures feasibility for any forecast error.  
Distributionally robust approaches use data to estimate distribution parameters (e.g., mean and variance) and aim to be robust to any data-generating distribution consistent with these parameters \cite{tyler,roald3,LiMathieu,Baker,DallAnese1,xie2018distributionally}. Others take a robust approach, assuming only knowledge of bounds on forecast errors and enforcing constraints for any possible realization, e.g., \cite{warrington,jabr}. Overall, this line of recent research has explored tractable approximations and reformulations of difficult stochastic optimal power flow problems. However, none of the existing work explicitly accounts for sampling errors arising from limited data, which in operation can cause poor out-of-sample performance\footnote{Out-of-sample performance is an evaluation of the optimal decisions using a dataset that is different from the one used to obtain the decision, which can be tested with Monte Carlo simulation}. Even with sophisticated recent stochastic programming techniques, decisions can be overly dependent on small amounts of relevant data from a high-dimensional space, a phenomenon akin to overfitting in statistical models.

We propose a multi-period data-based method to solve a stochastic optimal power flow problem based on limited information about forecast error distributions available through finite historical training datasets. A preliminary version of this work appeared in \cite{guo2017stochastic}, and here we significantly expand the work in several directions into the present two part paper. The main contributions are as follows
\vspace{-0.005cm}
\begin{enumerate}
\item We formulate a multi-stage distributionally robust optimal control problem for optimal power flow. A distributionally robust model predictive control algorithm is then proposed, which utilizes computationally tractable data-driven distributionally robust optimization techniques \cite{mohajerin} to solve the subproblems at each stage. Whereas distributionally robust optimization approaches focus on single-stage problems, here we extend these approaches to multi-stage settings to obtain closed-loop feedback control policies. This allow us to update forecast error datasets, and in turn re-compute decisions with the latest knowledge. In principle, the framework allows any forecasting methodology and a variety of ambiguity set parameterizations. We focus on Wasserstein balls \cite{givens} around an empirical data-based distribution \cite{givens,mohajerin}, which allows controllable conservativeness by adjusting the Wasserstein radius. In contrast to previous work, we obtain policies that are explicitly robust to sampling errors inherent in the dataset. This approach achieves superior out-of-sample performance guarantees in comparison to other stochastic optimization approaches, effectively regularizing against overfitting the decisions to limited available data.

\item We leverage pertinent linear approximations of the AC power-flow equations (see, e.g.,~\cite{Baran89,bolognani2015fast,guggilam2016scalable,christ2013sens,linModels}) to facilitate the development of computationally-affordable chance-constrained AC OPF solutions that are robust to distribution mismatches, and provide a unified framework that is applicable to both transmission and distribution systems. Formulations for distribution networks  incorporate inverted-based RESs and  energy storage systems, and focus on addressing the voltage regulation problem under uncertainty. The transmission system formulation incorporates synchronous generators and  power injections from RESs, and it focuses on probabilistic $N-1$ security constraints on active transmission line flows. The framework yields set points and feedback control policies for controllable devices that are robust to variations in solar and wind injections and sampling errors inherent to the finite training datasets.

\item The effectiveness and flexibility of the proposed methodologies are demonstrated with extensive numerical experiments in a 37-node distribution network and in a 118-bus transmission system. These extensive case studies are presented in Part II \cite{Guo2018DataDriven2}. We demonstrate inherent tradeoffs between economic efficiency and robustness to constraint violations and sampling errors due to forecasting. By explicitly incorporating forecast error and sampling uncertainties, the methodology can help network operators to better understand these risks and inherent tradeoffs, and to design effective optimization and control strategies for appropriately balancing efficiency objectives with security requirements.

\end{enumerate}

The rest of the paper is organized as follows: Section II presents a general formulation of the proposed data-based distributionally robust stochastic OPF problem. Section III describes the modeling of network, grid-connected components, and network constraints.
Section IV specializes the data-based stochastic OPF for distribution networks with an approximate AC power flow. Section V adapts the proposed methods for transmission systems with DC power flow. Section VI concludes the paper. Extensive simulation results from numerical experiments and discussions are presented in Part II.

\textbf{Notation}: The inner product of two vectors $a,b \in \mathbf{R}^m$  is denoted by $\langle a, b \rangle := a^\intercal b$. The $N_s$-fold product of distribution $\mathds{P}$ on a set $\Xi$ is denoted by $\mathds{P}^{N_s}$, which represents a distribution on the Cartesian product space $\Xi^{N_s} = \Xi \times \ldots \times \Xi$. We use $N_s$ to represent the number of samples inside the training dataset $\hat{\Xi}$. Superscript `` $\hat{\cdot}$ " is reserved for the objects that depend on a training dataset $\hat{\Xi}^{N_s}$. We use $(\cdot)^\intercal$ to denote  vector or matrix transpose.
The operators $\Re\{\cdot\}$ and $\Im\{\cdot\}$ return the real and imaginary part of a complex number, respectively. The operator $[\,\cdot\,]_{[a,b]}$ selects the $a$-th to $b$-th elements of a vector or rows of a matrix.
\section{Stochastic OPF as stochastic optimal control}
In this section, we formulate a stochastic OPF problem as a distributionally robust stochastic optimal control problem. We first pose the problem generically to highlight the overall approach, and in subsequent sections we detail the model and objective and constraint functions more explicitly for both distribution networks and transmission systems. This framework is more general than most stochastic OPF and distributionally robust optimization approaches in the literature, which typically focus only on individual or single-stage optimization problems.

Let $x_t \in \mathbf{R}^n$ denote a state vector at time $t$ that includes the internal states of all devices in the network. Let $u_t \in \mathbf{R}^m$ denote a control input vector that includes inputs for all controllable devices in the network. Let $\xi_t \in \mathbf{R}^{N_\xi}$ denote a random vector in a probability space $(\Omega, \mathcal{F}, \mathds{P}_t)$ that includes forecast errors of all uncertainties in the network. Forecast error distributions are never known in practice, so $\mathds{P}_t$ is assumed to be \emph{unknown} but belonging to an \emph{ambiguity set} $\mathcal{P}_t$ of distributions with a known parameterization, which will be discussed in detail shortly. We define the concatenated forecast error over an operating horizon $T$ as $\bm{\xi}_{0:T} := [\xi_0^\intercal,\ldots,\xi_{T}^\intercal]^\intercal \in \mathbf{R}^{N_\xi T}$, which has joint distribution $\mathds{P}$ and corresponding ambiguity set $\mathcal{P}$.

Since forecast errors are explicitly included, we seek closed-loop feedback policies of the form $u_t = \pi(x_0, \ldots, x_t, \bm{\xi}_{0:t}, \mathcal{D}_t)$, where the term $\mathcal{D}_t$ collects all network and device model information and the parameterization of the ambiguity set of the forecast error distribution. The control decisions at time $t$ are allowed to be functions of the entire state and forecast error history up to time $t$; this is called a history-dependent state and disturbance feedback information pattern. This general formulation allows for the design of not only nominal plans for controllable devices inputs, but also for planned reactions to forecast errors as they are realized.\footnote{The reactions can be interpreted as pre-planned secondary frequency control allocations \cite{warrington} or contingency reactions in response to forecast errors based on device dynamics and parameters describing forecast error distributions.} The policy function $\pi$ maps all available information to control actions and is an element of a set $\Pi$ of measurable functions.

This leads us to the following multi-stage distributionally robust stochastic optimal control problem
\begin{subequations} \label{generalstochasticopf}
\begin{align}
	&\inf_{\pi \in \Pi} \sup_{\mathds P \in \mathcal P} && \mathds{E}^{\mathds{P}} \sum_{t = 0}^{T} h_t(x_t, u_t, \xi_t), \\
	& \text{subject to} && x_{t + 1}  = f_t (x_t, u_t, \xi_t),\\
	& &&u_t = ~ \pi(x_0, \ldots, x_t, \bm{\xi}_t, \mathcal{D}_t),\\
	& && (x_t, u_t) \in \mathcal{X}_t.
\end{align}
\end{subequations}
The goal is to compute a feedback policy that minimizes the expected value of an objective function $h_t: \mathbf{R}^n \times \mathbf{R}^m \times \mathbf{R}^{N_\xi} \rightarrow \mathbf{R}$ under the worst-case distribution in the forecast error ambiguity set $\mathcal{P}$. The objective function $h_t$ will include both operating costs and risks of violating various network and device constraints and is assumed to be continuous and convex for every fixed $\xi_t$. The system dynamics function  $f_t: \mathbf{R}^n \times \mathbf{R}^m \times \mathbf{R}^{N_\xi t} \to \mathbf{R}^n$ models internal dynamics and other temporal interdependencies of devices, such as state of charge for batteries and ramp limits of generators. The constraint set $\mathcal{X}_t$ includes network and device constraints, such as power balance and generator and storage device bounds (some constraints may be modeled deterministically and others may be included as risk terms in the objective function).

The main challenges to solving \eqref{generalstochasticopf} are the multi-stage feedback nature of the problem, the infinite dimensionality of the control policies, the possible nonlinearity of device dynamics, and how to appropriately parameterize and utilize our available knowledge about forecast error distributions. We will tackle these using a distributionally robust model predictive control scheme with affine feedback policies and linear models for device dynamics, where stage-wise distributionally robust stochastic optimization problems are repeatedly solved over a planning horizon. To incorporate forecast error knowledge, we will use tractable reformulations of ambiguity sets based explicitly on empirical training datasets of forecast errors.


\subsection{Ambiguity sets and Wasserstein balls}
There is a variety of ways to reformulate the general stochastic OPF problem \eqref{generalstochasticopf} to obtain tractable subproblems that can be solved by standard convex optimization solvers. These include assuming specific functional forms for the forecast error distribution (e.g., Gaussian) and using specific constraint risk functions, such as those encoding value at risk (i.e., chance constraints), conditional value at risk (CVaR), distributional robustness, and support robustness. In all cases, the out-of-sample performance of the resulting decisions in operational practice ultimately relies on  1) the quality of data describing the forecast errors and 2) the validity of assumptions made about probability distributions. Many existing approaches make either too strong or too weak assumptions that possible lead to underestimation or overestimation of risk. In this paper, we extend a recently proposed tractable method \cite{mohajerin} to a multi-period data-based stochastic OPF, in which the ambiguity set is based on a finite forecast error training dataset $\hat{\Xi}^{N_s}$.


Within the  area of distributionally-robust optimization, moment-based ambiguity sets are utilized to model distributions featuring specified moment constraints such as unimodality \cite{LiMathieu,li2017ambiguous}, directional derivatives \cite{wiesemann2014distributionally}, symmetry \cite{roald3}, and log-concavity \cite{LiMathieu2018logconvae}. The ambiguity sets can be defined as confidence regions based on goodness-of-fit tests \cite{bertsimas2017robust}. Another line of works considers ambiguity sets as balls in the probability space, with radii computed based on the Wasserstein metric \cite{wozabal2012framework}, the Kullback-Leilber divergence \cite{jiang2016data}, and the Prohorov metric \cite{erdougan2006ambiguous}. This paper formulates ambiguity sets by leveraging Wasserstein balls. Relative to other approaches, Wasserstein balls provide an upper confidence bound, quantified by Wasserstein radius $\varepsilon$ \cite{mohajerin}, to achieve the superior out-of-sample performance; they also enable power system operators to ``control'' the conservativeness of the solution, thus ensuring he flexibility in the power system operation. Additionally, the approach in this paper seeks  the worst-case expectation subjected to all distributions contained in the ambiguity set. The worst-case expectation of the stochastic OPF problem over Wasserstein ambiguity set can be reformulated as a finite-dimensional convex problem, and can be solved using existing convex optimization solvers. 

The Wasserstein metric defines a distance in the space $\mathcal{M}(\Xi)$ of all probability distributions $\mathds{Q}$ supported on a set $\Xi$ with $\mathds{E}^{\mathds{Q}}[\| \xi \|] = \int_{\Xi} \| \xi \| \mathds{Q}(d\xi) < \infty$. In this paper, we assume the support set is polytopic of the form the uncertainty set is a polytope $\Xi := \{\xi \in \mathbf{R}^{N_\xi}: H\xi \leq \mathbf{d} \}$.

\textbf{Definition} [Wasserstein Metric]. Let $\mathcal{L}$ be the space of all Lipschitz continuous functions $f: \Xi \to \mathbf{R}$ with Lipschitz constant less than or equal to 1. The Wasserstein metric $d_W : \mathcal{M}(\Xi) \times \mathcal{M}(\Xi) \to \mathbf{R}$ is defined $\forall \mathds{Q}_{1}, \mathds{Q}_{2} \in \mathcal{M}(\Xi)$ as
\begin{equation}\label{WassersteinMetric}\nonumber
\begin{array}{l}
d_W(\mathds{Q}_{1}, \mathds{Q}_{2}) = \sup_{f \in \mathcal{L}}\bigg( \int_{\Xi} f(\xi)\mathds{Q}_{1}(d\xi) - \int_{\Xi} f(\xi)\mathds{Q}_{2}(d\xi) \bigg).
\end{array}
\end{equation}

Intuitively, the Wasserstein metric quantifies the minimum ``transportation'' cost to move mass from one distribution to another. We can now use the Wasserstein metric to define an ambiguity set
\begin{equation} \label{ambiguityset}
\hat{\mathcal{P}}^{N_s} := \bigg\{ \mathds{Q} \in \mathcal{M}(\Xi): d_w(\hat{\mathds{P}}^{N_s}, \mathds{Q}) \leq \varepsilon \bigg\},
\end{equation}
which contains all distributions within a Wasserstein ball of radius $\varepsilon$ centered at a uniform empirical distribution $\hat{\mathds{P}}^{N_s}$ on the training dataset $\hat{\Xi}^{N_s}$. The radius $\varepsilon$ can be chosen so that the ball contains the true distribution $\mathds{P}$ with a prescribed confidence level and leads to performance guarantees \cite{mohajerin}. The radius $\varepsilon$ also explicitly controls the conservativeness of the resulting decision. Large $\varepsilon$ will produce decisions that rely less on the specific features of the dataset $\hat \Xi^{N_s}$ and give better robustness to sampling errors. This parameterization will be used in the next subsection to formulate a distributionally robust optimization.
\subsection{Data-based distributionally robust model predictive control}
The goal of our data-based distributionally robust stochastic OPF framework is to understand and to illustrate inherent tradeoffs between efficiency and risk of constraint violations. Accordingly, the objective function comprises a weighted sum of an operational cost function and a constraint violation risk function: $h_t =  J_\textrm{Cost}^t + \rho J_\textrm{Risk}^t$, where $\rho \in \mathbf{R}_+$ is a weight that quantifies the network operator's risk aversion.
The operational cost function is assumed to be linear or convex quadratic
\begin{equation}\label{costfunction}
\nonumber J^t_{\text{Cost}}(x_t, u_t) := f_{x}^\intercal x_t +  \frac{1}{2}{x_t}^\intercal H_{x} x_t + f_{u}^\intercal u_t +  \frac{1}{2}{u_t}^\intercal H_{u}u_t,
\end{equation}
where $H_{x}$ and $H_{u}$ are positive semidefinite matrices. This function can capture several objectives including thermal generation costs, ramping  costs, and active power losses.

The risk function $J_\textrm{Risk}$ associated with the constraint violation  comprises a sum of the conditional value-at-risk (CVaR) \cite{rockafellar2000optimization} of a set of $N_\ell$ network and device constraint functions; specifically, it is defined as
\begin{equation}\label{generalDRO}\nonumber
J_\textrm{Risk}^t := \sum_{i=1}^{N_\ell} \text{CVaR}_{\mathds{P}}^\beta [\ell_i(x_t, u_t, \xi_t) ],
\end{equation}
where $\beta \in (0,1]$ refers to the confidence level of the CVaR under the distribution $\mathds{P}$ of the random variable $\xi_t$. Intuitively, the constraint violation risk function $J_{\textrm{Risk}}$ could be understood as the sum of networks and devices constraint violation magnitude at a ``risk level'' $\beta$, which penalizes both frequency and expected severity of constraint violations~\cite{rockafellar2000optimization}. Further details will be provided in subsequent sections. 

\textbf{Data-based distributionally robust model predictive control for stochastic OPF.} The general problem \eqref{generalstochasticopf} will be approached with a distributionally robust model predictive control (MPC) algorithm. MPC is a feedback control technique that solves a sequence of open-loop optimization problems over a planning horizon $\mathcal{H}_t$ (which in general may be smaller than the overall horizon $T$). At each time $t$, we solve the distributionally robust optimization problem over a set $\Pi_{\text{affine}}$ of affine feedback policies using the Wasserstein ambiguity set \eqref{ambiguityset}
\begin{subequations} 
\label{drmpcstochasticopf}
\begin{align}
	&\inf_{\pi \in \Pi_{\text{affine}}} \sup_{\mathds P \in \hat{\mathcal{P}}^{N_s} } && \mathds{E}^{\mathds{P}} \sum_{\tau = t}^{t+\mathcal{H}_t} J^\tau_{\text{Cost}}(x_\tau, u_\tau, \xi_\tau) + \rho J^\tau_{\text{Risk}}(x_\tau, u_\tau, \xi_\tau),  \\
	& \text{subject to} && x_{\tau + 1}  = f_\tau (x_\tau, u_\tau, \xi_\tau), \\
	& &&u_\tau = ~ \pi(x_0, \ldots, x_\tau, \bm{\xi}_\tau, \mathcal{D}_\tau), \\
	& && (x_\tau, u_\tau) \in \mathcal{X}_\tau.
\end{align}
\end{subequations}
Only the immediate control decisions for time $t$ are implemented on the controllable device inputs. Then time shifts forward one step, new forecast errors and states are realized, the optimization problem \eqref{drmpcstochasticopf} is re-solved at time $t+1$, and the process repeats. This approach allows any forecasting methodology to be utilized to predict uncertainties over the planning horizon. Furthermore, the forecast error dataset $\hat{\mathds{P}}^{N_s}$, which defines the center of the ambiguity set $\hat{\mathcal{P}}^{N_s}$, can be updated online as more forecast error data is obtained. It is also possible to remove outdated data online to account for time-varying distributions.

In the rest of the paper, we will derive specific models for both distribution and transmission networks and grid devices where the subproblems \eqref{drmpcstochasticopf} have exact tractable convex reformulations as quadratic programs \cite{mohajerin} and can be solved to global optimality with standard solvers.

\section{System model}
We now consider a symmetric and balanced electric power network model in steady state, where all currents and voltages are assumed to be sinusoidal signals at the same frequency.
\subsection{Network model}
Consider a power network (either transmission or distribution), denoted by a graph $(\mathcal{N},\mathcal{E})$, with a set $\mathcal{N} = \{1,2,...,N\}$ of buses, and a set $\mathcal{E} \subset \mathcal{N} \times \mathcal{N}$ of the power lines connecting buses. Let $V_i^t \in \mathbf{C}$ and $I_i^t \in \mathbf{C}$ denote the phasors for the line-to-ground voltage and the current injection at node $i \in \mathcal{N}$. Define the complex vectors $\mathbf{v}^t:=[V_1^t, V_2^t, ..., V_N^t]^\intercal \in \mathbf{C}^N$ and $\mathbf{i}^t :=[I_1^t,...,I_N^t] \in \mathbf{C}^N$. Let $z_{ij}$ denote the complex impedance of the line between bus $i$ and bus $j$, then the line admittance is $y_{ij} = 1/z_{ij} = g_{ij} + \textrm{j}b_{ij}$. We model the lines using a standard Pi Model. The admittance matrix $\mathbf{Y} \in \mathbf{C}^{N \times N}$ has elements
\begin{equation}\label{admittanceMatrix}
Y_{ij} = \left\{ \begin{array}{ll}
\sum_{l \sim i} y_{il} + y_{ii} & \textrm{if $i=j$}\\
-y_{ij} & (i,j) \in \mathcal{E} \\
0 &  (i,j) \notin \mathcal{E}
\end{array} \right .,
\end{equation}
where $l \sim i$ indicates that bus $i$ and bus $j$ are connected. Via Kirchoff's and Ohm's laws, we have $\mathbf{i}^t = \mathbf{Y}\mathbf{v}^t.$
Net complex power bus injections are given by
\begin{equation} \label{powerflow}
\mathbf{s}^t = \mathbf{v}^t \, (\mathbf{i}^t)^* = \textrm{diag}\left(\mathbf{v}^t
\right) \left(\mathbf{Y}\mathbf{v}^t\right)^*.
\end{equation}
The components of $\mathbf{s}^t = [S_1^t, S_2^t, \ldots , S_N^t]^\intercal \in \mathbf{C}^N$ can be expressed in rectangular coordinates as $S_i^t = P_i^t + \textrm{j}Q_i^t$, where $P_i^t$ is active power and $Q_i^t$ is reactive power. Positive $P_i$ and $Q_i$ means that bus $i$ generates active/reactive power, and negative $P_i$ and $Q_i$ mean that bus $i$ absorbs the active/reactive power. Vectors of active and reactive power $\mathbf{p}^t = [P_1^t , P_2^t, \ldots ,P_N^t]^\intercal$ and $\mathbf{q}^t = [Q_1^t, Q_2^t, \ldots ,Q_N^t]^\intercal$ are further divided into nominal and error terms: $\mathbf{p}^t = \bar{\mathbf{p}}^t(u_t) + \tilde{\mathbf{p}}^t(\xi_t)$ and $\mathbf{q}^t = \bar{\mathbf{q}}^t(u_t) + \tilde{\mathbf{q}}^t(\xi_t)$. The nominal active and reactive power injection vectors $\bar{\mathbf{p}}^t(u_t) \in \mathbf{R}^N$ and $\bar{\mathbf{q}}^t(u_t) \in \mathbf{R}^N$ depend on control decisions, and the forecast errors $\tilde{\mathbf{p}}^t(\xi_t) \in \mathbf{R}^N$ and $\tilde{\mathbf{q}}^t(\xi_t)  \in \mathbf{R}^N$ depend on the random vector $\xi_t$.

To handle nonconvexity of the power flow equations \eqref{powerflow}, we utilize two different linearization methods that are effective in both distribution and transmission networks.

\subsection{Dynamic of grid-connected devices}
We consider $N_d$ grid-connected devices, which may include 1) traditional generators and inverter-based RESs; 2) fixed, deferrable, and curtailable loads; 3) storage devices like batteries and plug-in electric vehicles, which are able to act as both generators and loads. There are two types of devices: devices with controllable power flow affected by decision variables (e.g., conventional thermal and RESs generators, deferrable/curtailable loads and storage devices); and devices  with fixed or uncertain power flow which will not be affected by decision variables (e.g., fixed loads).
The power flow of each controllable device is modeled by a discrete-time linear dynamical system
\begin{equation}\label{dynamic}\nonumber
x_{t+1}^d = \bar{A}^d x_{t}^d + \bar{B}^d u_{t}^d,
\end{equation}
where device $d$ at time $t$ has internal state $x_t^d \in \mathbf{R}^{n_d}$, dynamics matrix $\bar{A}^d \in \mathbf{R}^{n_d \times n_d}$, input matrix $\bar{B}^d \in \mathbf{R} ^{n_d \times m_d}$, and control input $u_t^d \in \mathbf{R}^{m_d}$. The first element of $x_t^d$ corresponds to the power injection of device $d$ at time $t$ into the network at its bus, and other elements describe internal dynamics, such as state-of-charge (SOC) of  energy storage devices. At time $t$, state and input histories are expressed by $\mathbf{x}^d_t : = [(x_1^d)^\intercal, \ldots , (x_t^d)^\intercal ]^\intercal \in \mathbf{R}^{n_d t} $ and $\mathbf{u}^d_t : = [(u_0^d)^\intercal, \ldots , (u_{t-1}^d)^\intercal ]^\intercal \in \mathbf{R}^{m_d t}$ with
\begin{equation}\label{statefunction}\nonumber
\mathbf{x}^d_t = A^d_t x_0^{d} + B^d_t\mathbf{u}^d_t,
\end{equation}
where
\begin{equation*}
A^d_t :=
\begin{bmatrix}
\bar{A}^d\\
(\bar{A}^d)^2\\
\vdots\\
(\bar{A}^d)^t
\end{bmatrix}
,
B^d_t :=
\begin{bmatrix}
\bar{B}^d & 0 & \dots & 0\\
\bar{A}^d\bar{B}^d & \bar{B}^d& \ddots & 0 \\
\vdots & \ddots & \ddots & \vdots \\
(\bar{A}^d)^{t-1} \bar{B}^d & \dots & \bar{A}^d\bar{B}^d & \bar{B}^d
\end{bmatrix}.
\end{equation*}
\subsection{Network constraints}
The AC power flow equations render prototypical AC OPF formulation nonconvex and NP-hard; what is more, in the present context, they hinder the development of computationally-affordable chance-constrained AC OPF formulations where CVaR arguments are leveraged as risk measures. 

For distribution systems, we refer the reader to the linear approximation methods proposed in e.g.,~\cite{Baran89,bolognani2015fast,guggilam2016scalable,christ2013sens,linModels,dhople2015linear}, with the latter suitable for multi-phase systems with both wye and delta connections; these approximate models have been shown to provide high levels of accuracy in many existing test systems.  For transmission systems, one can consider the tradition DC power-flow model to approximate the voltage angles and active power flows in the system; see, e.g.,~\cite{christie2000transmission,warrington,tyler,wood2012power}. Alternatively, one can consider alternative linearizations; see e.g.,\cite{dhople2015linear}. As long as an accurate linear model exists, the proposed technical approach can be utilized to formulate and solve a distributionally-robust chance-constrained AC OPF problem. In particular, linear models can be utilized to formulate convex (in fact, linear) constraints on \emph{line flows} and \emph{voltage magnitudes}, e.g.  
\begin{equation}\label{generalvoltagemagnitudeconstraints}
V^{\rm{min}} \leq |V_i| \leq V^{\rm{max}}, \,\, \forall i\in \mathcal{N}.
\end{equation}

Grid-connected devices have various local constraints including, e.g., state of charge limitations for energy storage devices, allowable power injection ranges, generator ramping limits,  and other device limits.  These can be modeled (or approximated) as linear inequalities of the form
\begin{equation}\label{localdeviceconstraints} 
\mathbf{T}_d^t \mathbf{x}^d_t + \mathbf{U}_d^t\mathbf{u}^d_t + \mathbf{Z}_d^t\bm{\xi}_t \leq w_d, \,\,\, d = 1, \ldots , N_d,
\end{equation}
where $\mathbf{T}_d^t \in \mathbf{R}^{l_d \times n_d t}$, $\mathbf{U}_d^t \in \mathbf{R}^{l_d \times m_d t}$, and $\mathbf{Z}_d^t \in \mathbf{R}^{l_d \times N_\xi t}$, and $w_d\in \mathbf{R}^{l_d}$ is a local constraint parameter vector.

\section{Data-based distributionally robust \\stochastic OPF for distribution networks}
\subsection{Distribution network model}
In this section, we specialize the model to symmetric and balanced power distribution networks, connected to the grid at a point of common coupling (PCC). Loads and distributed generators (e.g., thermal generators, inverter-based RESs, and energy storage devices) may be connected to each bus. We augment the bus set with node $0$ as the PCC.

The voltage and injected current at each bus are defined as $V^t_n = |V^t_n|e^{j\angle V^t_n}$, and $I^t_n = |I^t_n|e^{j\angle I^t_n}$. The absolute values $|V^t_n|$ and $|I^t_n|$ correspond to the signal root-mean-square values, and the phase $\angle V_n^t$ and phase $\angle I_n^t$ correspond to the phase of the signal with respect to a global reference.

Node 0 is modeled as a slack bus and the others are PQ buses, in which the injected complex power are specified. The admittance matrix can be partitioned as
\begin{equation}\nonumber
\begin{bmatrix}
I_0^t\\
\mathbf{i}^t
\end{bmatrix} =
\begin{bmatrix}
y_{00} & \bar{y}^\intercal  \\
\bar{y} & \mathbf{Y}
\end{bmatrix}
\begin{bmatrix}
V_0\\
\mathbf{v}^t
\end{bmatrix}.
\end{equation}

The net complex power injection is then
\begin{equation}\label{powerbalancewithPCC}
\mathbf{s}^t = \textrm{diag}(\mathbf{v}^t)\Big(\mathbf{Y}^*(\mathbf{v}^t)^* + \bar{y}^*(v_0^t)^* \Big).
\end{equation}
The nonconvexity of this equation in the space of power injections and bus voltages is a source of significant computational difficulty in optimal power flow problems. In the rest of this section, we formulate a convex and computationally efficient data-based stochastic OPF problem based on a particular linear approximation of \eqref{powerbalancewithPCC} that is appropriate for distributions networks. This approximation occurs on a specific point of a power flow manifold that preserves network structure for both real and reactive power flows and allows direct application of stochastic optimization techniques for incorporating forecast errors.


\subsection{Leveraging approximate power flow} 

Collect the voltage magnitudes $\{|V_n^t| \}_{n\in \mathcal{N}}$ into the vector $|\mathbf{v}^t| := [|V_1^t|,
\ldots, |V_N^t|]^\intercal \in \mathbf{R}^N$. To develop computationally-feasible approaches, the technical approach in this paper leverages linear approximations of the AC power-flow equations; in particular, linear approximations for voltages, as a function of the injected power $\mathbf{s}^t$, are given by
\begin{equation}\label{rectangularpolar}
\mathbf{v}^t \thickapprox \bigg(\mathbf{H}\mathbf{p}^t + \mathbf{J}\mathbf{q}^t + \mathbf{c}\bigg), \quad |\mathbf{v}^t| \thickapprox \bigg(\mathbf{M}\mathbf{p}^t + \mathbf{N}\mathbf{q}^t + \mathbf{a}\bigg) . 
\end{equation}
Using these approximations, the voltage constraints $V^{\rm{min}} \leq |V_n^t| \leq V^{\rm{max}}$ can be approximated as  $V^{\rm{min}}\mathbf{1}_{N} \preceq \mathbf{M}\mathbf{p}^t + \mathbf{N}\mathbf{q}^t + \mathbf{a} \preceq V^{\rm{max}}\mathbf{1}_N.$ 
The coefficient matrices of the linearized voltages, and the normalized vectors $\mathbf{a}$ and $\mathbf{c}$ can be obtained as shown in \cite{dhople2015linear,guggilam2016scalable,linModels, christ2013sens}. For completeness, in the remainder of this subsection we briefly outline the approach taken in~\cite{dhople2015linear,guggilam2016scalable} to derive a linear model for the voltages. 

Suppose that $\mathbf{v}^t = \bar{\mathbf{v}} + \Delta v^t,$ where $\bar{\mathbf{v}} = |\bar{\mathbf{v}}| \angle \bm \theta  \in \mathbf{C}^N$ is a pre-determined nominal voltage vector and $\Delta v^t \in \mathbf{C}^N$ denotes a deviation from nominal. Then we have
\begin{equation}
\mathbf{s}^t = \textrm{diag}(\bar{\mathbf{v}} + \Delta v^t)\bigg(\mathbf{Y}^*(\bar{\mathbf{v}}+\Delta v^t)^* + \bar{y}^*V_0^*\bigg).
\end{equation}
Neglecting second-order terms $\textrm{diag}(\Delta v^t)\mathbf{Y}^* (\Delta v^t)^*$, the power balance \eqref{powerbalancewithPCC} becomes $\Lambda \Delta v^t + \Phi (\Delta v^t)^* = \mathbf{s}^t + \Psi$, where $\Lambda := \textrm{diag}(\mathbf{Y}^*\bar{\mathbf{v}}^*+\bar{y}^*V_0^*)$, $\Phi := \textrm{diag}(\bar{\mathbf{v}})\mathbf{Y}^* $, $\Psi := -  \textrm{diag}(\bar{\mathbf{v}})(\mathbf{Y}^*\bar{\mathbf{v}}^* + \bar{y}^*V_0^*)$.
We consider a choice of the nominal voltage
$\bar{\mathbf{v}} =  \mathbf{Y}^{-1}\bar{y}V_0,$ which gives $\Lambda = \mathbf{0}_{N \times N}$ and $\Psi = \mathbf{0}_N$. Therefore the linearized power flow expression is
$\mathbf{s}^t = \textrm{diag}(\bar{\mathbf{v}})\mathbf{Y}^*(\Delta v^t)^*$,
the deviation $\Delta v^t$ becomes
$\Delta v^t = \mathbf{Y}^{-1}\textrm{diag}^{-1}(\bar{\mathbf{v}}^*)(s^t)^*.$

Let us denote $\mathbf{Y}^{-1} = (\mathbf{G} + \textrm{j}\mathbf{B})^{-1} = \mathbf{Z}_R + \textrm{j}\mathbf{Z}_I$. Then expanding $\Delta v^t$ in rectangular form gives
\begin{equation*}
\begin{split}
\mathbf{\bar{M}} & = \bigg( \mathbf{Z}_R~ \textrm{diag}\bigg(\frac{\textrm{cos}(\theta)}{|\bar{\mathbf{v}}|}\bigg) - \mathbf{Z}_I~ \textrm{diag}\bigg(\frac{\textrm{sin}(\theta)}{|\bar{\mathbf{v}}|}\bigg)\bigg),\\
\mathbf{\bar{N}} & = \bigg( \mathbf{Z}_I~ \textrm{diag}\bigg(\frac{\textrm{cos}(\theta)}{|\bar{\mathbf{v}}|}\bigg) + \mathbf{Z}_R~ \textrm{diag}\bigg(\frac{\textrm{sin}(\theta)}{|\bar{\mathbf{v}}|}\bigg)\bigg),\\
\end{split}
\end{equation*}
which define the rectangular matrices $\mathbf{H} := \mathbf{\bar{M}} + j \mathbf{\bar{N}}$, $\mathbf{J} := \mathbf{\bar{N}} - j \mathbf{\bar{M}}$, and the coefficient $\mathbf{c}$ is $\bar{\mathbf{v}}$. If $\bar{\mathbf{v}}$ dominates $\Delta v^t $, then the voltage magnitudes are approximated by $|\bar{\mathbf{v}}| + \Re \{\Delta v^t\}$, and linearized coefficients of voltage magnitudes become $\mathbf{M} :=\mathbf{\bar{M}}$, $\mathbf{N} := \mathbf{\bar{N}}$, and $\mathbf{a} := |\mathbf{\bar{v}}|$. 
It is worth noting that the approach proposed in~\cite{linModels} accounts for multiphase systems with both wye and delta connections. Accordingly, the proposed framework is applicable to generic multiphase feeders with both wye and delta connections.

\subsection{Data-based stochastic OPF formulation for distribution networks}
Using the introduced linearized relationship between voltage and power injection vectors $\mathbf{p}^t$ and $\mathbf{q}^t$, we express the voltage magnitude in the following form
\begin{equation} \nonumber \label{linearvoltageapproximation}
g^t\big[\mathbf{p}^t(u_t,\xi_t), \mathbf{q}^t(u_t, \xi_t) \big]: = \mathbf{M}(\mathbf{I}-\textrm{diag}\{\bm{\alpha}^t\})\mathbf{p}^t_{\rm{av}} + \mathbf{N}\mathbf{q}^t + |\bar{\mathbf{v}}|,
\end{equation}
where $\bm{\alpha}^t \in \mathbf{R}^N$ is a control policy defined as the fraction of the active power curtailment by the renewable energy power injection. A system state vector  $\mathbf{p}_{\rm{av}}^t \in \mathbf{R}^N$ collects the active power injection including loads and the available RES power. We aim to optimize the set points $\{\bm{\alpha^t}, \mathbf{q}^t \}$ of nodal power injections, which can be achieved by adjusting controllable loads and generators. More details of system modeling and component dynamics will be introduced in Part II. 

Broadly speaking, we quantify a violation risk of voltage magnitude constraints (6) and local device constraints (7) for each node and each time as follows
\begin{subequations}
\begin{align}
& \hspace{-3mm} \mathds{E}~\mathcal{R}\bigg\{g_{n}^t \big[\mathbf{p}^t(u_t,\xi_t), \mathbf{q}^t(u_t, \xi_t) \big] - V^{\rm{max}} \bigg\} \leq 0,~ \label{PRVoltageMax} \\
&\hspace{-3mm} \mathds{E}~\mathcal{R} \bigg\{V^{\rm{min}} - g_{n}^t\big[\mathbf{p}^t(u_t,\xi_t), \mathbf{q}^t(u_t, \xi_t) \big]\bigg\}  \leq 0,\label{PRVoltageMin} \\
& \hspace{-3mm} \mathds{E}~\mathcal{R} \bigg\{\mathbf{T}_d^t \mathbf{x}_t^d + \mathbf{U}_d^t \mathbf{u}^d_t + \mathbf{Z}_d^t \bm{\xi}_t - w_d \bigg\} \leq 0 , d = 1,\ldots, N_d, \label{PRLocalDevices}
\end{align}
\end{subequations}
where $g_n^t{(\,\cdot\,)}$ is the $n$-th element of the function value $g^t(\,\cdot \,)$, and $\mathcal{R}$ denotes a generic transformation of the inequality constraints into stochastic versions. Using a prior on the uncertainty and possibly introducing auxiliary variables, the general risk functions \eqref{PRVoltageMax}-\eqref{PRLocalDevices} can be reformulated using e.g., scenario-based approaches,  \cite{vrakopoulou3,li2017chance} or moment-based distributionally robust optimization \cite{tyler,LiMathieu}. This paper seeks a reformulation by leveraging the CVaR~\cite{rockafellar2000optimization}. A set of constraints will be approximated using the proposed distributionally robust approach, while  other constraints will be evaluated using sample average methods. 

\color{black}
We define a set $\mathcal{V}_t$ that contains a subset of $N_\ell$ affine constraints \eqref{PRVoltageMax}-\eqref{PRLocalDevices} that will be treated with distributionally robust optimization techniques. This allows some or all of the constraints to be included. We express them in terms of a decision variable vector $\mathbf{y}_t$ and uncertain parameters $\xi_t$, where $\mathbf{y}_t$ consists of all the decision variables including the RES curtailment ratio vector $\bm{\alpha}^\tau$ and other controllable device set-points, and $\xi_t$ contains the uncertain parameters across the network including the active and reactive power injection forecast errors. For simplicity, we consider the risk of each constraint individually; it is possible to consider risk of joint constraint violations, but this is more difficult and we leave it for future work. Each individual affine constraint in the set $\mathcal{V}_t$ can be written in a compact form as follows
\begin{equation*}
\begin{split}
\mathcal{C}_o^t(\mathbf{y}_t,\xi_t) = [\bar{\mathcal{A}}(\mathbf{y}_t)]_o\xi_t + [\bar{\mathcal{B}}(\mathbf{y}_t)]_o, \quad o = 1,..., N_\ell,
\end{split}
\end{equation*}
where $\mathcal{C}_o^t(\cdot)$ is the $o$-th affine constraint in the set $\mathcal{V}_t$. We use $[\, \cdot \,]_o$ to denote the $o$-th element of a vector or $o$-th row of a matrix. The CVaR with risk level $\beta$ of the each individual constraint in the set $\mathcal{V}_t$ is
\begin{equation} \label{defininationofCVaRdistribution}
\inf_{\kappa_o^t} \mathds{E}_{\xi_t}\bigg\{[\mathcal{C}^t_o(\mathbf{y}_t,\xi_t) + \kappa_o^t]_+ - \kappa_o^t\beta \bigg\} \leq 0,
\end{equation}
where $\kappa_o^t \in \mathbf{R}$ is an auxiliary variable \cite{rockafellar2000optimization}. The expression inside the expectation in \eqref{defininationofCVaRdistribution} can be expressed in the form
\begin{equation}\label{distributionCVaRConvex}
\bar{\mathcal{Q}}_o^t = \max_{k=1,2} \bigg[ \langle \bar{\mathbf{a}}_{ok}(\mathbf{y}_t), \xi_t \rangle + \bar{\mathbf{b}}_{ok}(\kappa_o^t) \bigg].
\end{equation}
This expression is convex in $\mathbf{y}_t$ for each fixed $\xi_t$ since it is the maximum of two affine functions. Our risk objective function is expressed by the distributionally robust optimization of CVaR
\begin{equation}\nonumber
\begin{split}
& \hat{J}_\textrm{Risk}^t = \sum_{\tau = t}^{t+\mathcal{H}_t}\sum_{o=1}^{N_\ell} \sup_{\mathds{Q}_\tau \in \hat{\mathcal{P}}_\tau^{N_s}}\mathds{E}^{\mathds{Q}_\tau} \max_{k=1,2} \bigg[ \langle \bar{\mathbf{a}}_{ok} (\mathbf{y}_\tau), \hat{\xi}_\tau \rangle + \bar{\mathbf{b}}_{ok}(\kappa_o^\tau) \bigg].
\end{split}
\end{equation}
The above multi-period distributionally robust optimization can be equivalently reformulated the following quadratic program, the details of which are described in \cite{mohajerin}.
The objective is to minimize a weighted sum of an operational cost function and the total worst-case CVaR of the affine constraints in set $\mathcal{V}^t$ (e.g., voltage magnitude and local device constraints).\\ 
\hspace{-2mm}\textbf{Data-based distributionally robust stochastic OPF}
\begin{subequations} \label{DROOPFdistributionsystem}
\begin{align}
	&\hspace{-2mm}\inf_{\mathbf{y}_\tau, \kappa_o^\tau } \sum_{\tau = t}^{t + \mathcal{H}_t}  \bigg\{\mathds{E} [\hat{J}^\tau_{\text{Cost}}] +  \rho \hspace{-2 mm} \sup_{\mathds{Q}_\tau \in \hat{\mathcal{P}}^{N_s}_\tau} \sum_{o=1}^{N_\ell} \mathds{E}^{\mathds{Q}_\tau}[\bar{\mathcal{Q}}_o^\tau]  \bigg\},  \nonumber\\
	&\hspace{-2mm} = \inf_{\begin{subarray}{c}\mathbf{y}_\tau, \kappa_o^\tau, \\  \lambda_o^\tau, s_{io}^\tau, \varsigma_{iko}^\tau \end{subarray}} \sum_{\tau=t}^{t + \mathcal{H}_t} \bigg\{\mathds{E}[\hat{J}^t_{\text{Cost}}] + \sum_{o=1}^{N_\ell} \bigg(\lambda_o\varepsilon_\tau + \frac{1}{N_s}\sum_{i=1}^{N_s} s_{io}^\tau\bigg) \bigg\},\\
	&\hspace{-2mm} \text{subject to}~~~~~~~~~~~~~~~~~~~~~~~~~~~~~~~~~~~~~~~~~~~~~~~~~~~~~~~~~~~~~~~~~~~~~~~~~~\nonumber \\
	&\hspace{-2mm}\rho (\bar{\mathbf{b}}_{ok}(\kappa_o^\tau) + \langle \bar{\mathbf{a}}_{ok}(\mathbf{y}_\tau), \hat{\xi}_\tau^{i}\rangle + \langle\varsigma_{iko}, \mathbf{d}-H\hat{\xi}_\tau^i\rangle) \le s_{io}^\tau,\\
	&\hspace{-2mm} \|H^\intercal \varsigma_{iko}-\rho \bar{\mathbf{a}}_{ok}(\mathbf{y}_\tau)\|_\infty \le \lambda_o^\tau,\\
	&\hspace{-2mm}\varsigma_{iko} \geq 0,\\
	&\hspace{-2mm}\forall i\le N_s, \forall o \le N_\ell, k=1,2, \tau = t,..., t + \mathcal{H}_t, \nonumber 
\end{align}
\end{subequations}
where $\rho \in \mathbf{R}_+$ quantifies the power system operators' risk aversion. This is a quadratic program that explicitly uses the training dataset $\hat{\Xi}^{N_s}_\tau = \{\hat{\xi}_\tau^i\}_{i \leq N_s}$. The risk aversion parameter $\rho$ and the Wasserstein radius $\varepsilon_\tau$ allow us to explicitly balance tradeoffs between efficiency, risk and sampling errors inherent in $\hat{\Xi}^{N_s}_\tau$. The uncertainty set is modeled as a polytope $\Xi := \{\xi \in \mathbf{R}^{N_\xi}: H\xi \leq \mathbf{d} \}$. The constraint $\varsigma_{iko} > 0$ holds since the uncertainty set is not-empty; on the other hand, in a case with no uncertainty (i.e, $\varsigma_{iko} = 0$), the variable $\lambda$  does not play any role and $s_{io}^\tau = \rho (\bar{\mathbf{b}}_{ok}(\kappa_o^\tau) + \langle \bar{\mathbf{a}}_{ok}(\mathbf{y}_\tau), \hat{\xi}_\tau^{i}\rangle)$. 

\section{Data-based distributionally robust stochastic optimal power flow for transmission systems}
\subsection{Illustrative explanation for the DC approximation}
In transmission networks, we employ a widely used ``DC'' linearization of the nonlinear AC power flow equations that makes the following assumptions \cite{warrington,tyler}:
\begin{itemize}
\item The lines $(i,j) \in \mathcal{E}$ are lossless and characterized by their reactance $\Im\{1/y_{ik}\}$;
\item The voltage magnitudes $\{ |V_n| \}_{\{n \in \mathcal{N}\}}$ are constant and close to one per unit;
\item Reactive power is ignored.
\end{itemize}
Under a DC model, line flows are linearly related  to the nodal power injections; therefore,  flow constraints in the transmission lines can be expressed as
\begin{equation} 
\label{lineconstraints}
\sum_{d=1}^{N_d} \Gamma_d^{t} (r_d^{t} + G_d^{t}\bm{\xi}_{t} + C_d^{t}\mathbf{x}^d_{t}) \le \bar{p}_{t},
\end{equation}
where $\Gamma_d^{t} \in \mathbf{R}^{{2Lt}\times t}$ maps the power injection (or consumptions) et each node to line flows and  can be constructed from the network line impedances \cite{christie2000transmission,zimmerman2010matpower}; on the other hand,  $\bar{p}_{t} \in \mathbf{R}^{2Lt}$ denotes a limit on the line flows. The power injections in~\eqref{lineconstraints} features controllable and non-controllable components;  the  non-controllable power for a device/node $d$ is modeled as $r_d^{t} + G_d^{t}\bm{\xi}_{t}$ (with positive values denoting the net power injection in the network), where the vector $r_d^{t} \in \mathbf{R}^{t}$ denotes the nominal power injection and $G_d^{t} \bm{\xi}_{t}$ models uncertain injections. Finally, $\mathbf{x}_t^d$ is the vector of controllable powers, and the matrix $C_d^{t} \in \mathbf{R}^{t \times n_d t}$ selects the first element of $\mathbf{x}^d_{t}$ at each time, i.e., $C_d^t := I_t \otimes \tilde{C}_d^t$, where $\tilde{C}_d^t = [1, \mathbf{0}_{1\times (n_d -1)}]$,  $I_t$ is a $t$-dimensional identity matrix and $\otimes$ denotes the Kronecker product operator. The power balance constraint is
\begin{equation}\label{powerbalance}
\sum_{d=1}^{N_d} (r_d^{t} + G_d^{t}\bm{\xi}_{t} +C_d^{t}\mathbf{x}^d_{t}) = 0 . 
\end{equation}

\subsection{Reserve policies}
Deterministic OPF formulations ignore the uncertainties $\bm{\xi}_t$ and compute an open-loop input sequence for each device. In a stochastic setting, one must optimize over causal \emph{policies}, which specifies how each device should respond to the current system states and forecast errors as they are discovered. We can now formulate a finite horizon stochastic optimization problem

\begin{subequations}
\label{stochopf}
\begin{align}
	&\hspace{-2.5mm} \inf_{\pi \in \Pi_{\textrm{affine}}}  \sum_{\tau = t}^{t+\mathcal{H}_t} \mathds{E}\Big[J^\tau_{\textrm{Cost}}(x_\tau, \pi_\tau,\xi_\tau)\Big],\\
	&\hspace{-2.5mm} \textrm {subject to} \nonumber \\
	&\hspace{-2.5mm} \sum_{d=1}^{N_d} \Big[r_d^{\bar{t}} + G_d^{\bar{t}}\bm{\xi}_{\bar{t}} + C_d^{\bar{t}}\mathbf{x}_{\bar{t}}^d\Big]_{[\underline{t},\bar{t}\,]} = 0,\\
    &\hspace{-2.5mm}\mathds{E}~\mathcal{R}\bigg\{\sum_{d=1}^{N_d} \Big[\Gamma_d^{\bar{t}}(r_d^{\bar{t}} + G_d^{\bar{t}}\bm{\xi}_{\bar{t}} + C_d^{\bar{t}}\mathbf{x}_{\bar{t}}^d) -\bar{p}_{\bar{t}}\Big]_{[2L\underline{t},2L\bar{t}\,]}  \bigg\} \le 0, \label{TransmissionLineFlow}\\ 
    &\hspace{-2.5mm} \mathds{E}~\mathcal{R}\Big[\mathbf{T}_d^\tau \mathbf{x}_{\tau}^d + \mathbf{U}_d^\tau \mathbf{u}^d_{\tau}+ \mathbf{Z}_d^\tau \bm{\xi}_{\tau} - w_d \Big] \le 0,\label{TransmissionLocalDevices}\\
    &\hspace{-2.5mm} d = 1,\ldots, N_d, \tau = t, \ldots, t + \mathcal{H}_t, \nonumber
\end{align}
\end{subequations}
where $\mathcal{R}$ denotes a general constraint risk function. The definition and discussion of the general stochastic transformation of \eqref{TransmissionLineFlow}-\eqref{TransmissionLocalDevices} is same as \eqref{PRVoltageMax}-\eqref{PRLocalDevices}, and can be found in Section IV.C. Here, $[\underline{t}, \bar{t}]$ denotes the finite time interval $[t, t+\mathcal{H}_t]$ for brevity. The cost function is proportional to the first and second moments of the uncertainties $\xi_t$, since we assume the operational cost function is convex quadratic. Optimizing over general policies makes the problem \eqref{stochopf} infinite dimensional, so we restrict attention to affine policies
\begin{equation} \label{inputpolicy}\nonumber
\mathbf{u}^d_t = D^d_t \bm{\xi}_t + e^d_t,
\end{equation}
where each participant device $d$ (e.g., traditional generators, flexible loads or energy storage devices) power schedule $\mathbf{u}^d_t \in \mathbf{R}^{t}$ is parameterized by a nominal schedule $e^d_t \in \mathbf{R}^t$ plus a linear function $D_t^d \in \mathbf{R}^{t \times N_\xi t}$ of prediction error realizations. To obtain causal policies, $D_t^d$ must be lower-triangular. The $D_t^d$ matrices can be interpreted as planned reserve mechanisms involving secondary frequency feedback controllers, which adjust conventional generator power outputs using frequency deviations to follow power mismatches \cite{warrington}. Under affine policies, the power balance constraints are linear functions of the distribution of $\bm{\xi}_t$, which are equivalent to
\begin{equation*}
\sum_{d=1}^{N_d}(r_d^t + C_d^t(A_t^dx_0^d + B_t^de_t^d)) = 0, \sum_{d=1}^{N_d}(G_d^t + C_d^tB_t^dD_t^d) = 0.
\end{equation*}

\subsection{Data-based stochastic OPF formulation for transmission systems}
We now use the above developments to formulate a data-based distributionally robust OPF problem to balance an operating efficiency metric with CVaR values of line flow and local device constraint violations. We collect the line flow and local device constraints into a set $\mathcal{V}_t$, which will be evaluated with distributionally robust optimization techniques in a total amount of $N_\ell$. This allows some or all of the constraints to be included. For simplicity, we again consider the risk of each constraint individually. We express the constraints in terms of decision variables $\{ D_t^d, e_t^d\}$ and uncertain parameters $\xi_t$. Each individual affine constraint in the set $\mathcal{V}_t$ can be written in a compact form as follows
\begin{equation*} \label{linearconstraintgeneral}
\mathcal{C}_o^t(D_t,e_t,\xi_t) = [\underline{\mathcal{A}}(D_t)]_o\xi_t + [\underline{\mathcal{B}}(e_t)]_o, ~ o = 1,...,N_\ell,
\end{equation*}
where $\mathcal{C}_o^t( \cdot )$ is the $o$-th affine constraint in the set $V_t$. The decision variables $D_t$ and $e_t$  both appear linearly in $\mathcal{C}_o^t$. The CVaR with risk level $\beta$ of the each individual constraint in the set $\mathcal{V}_t$ is
\begin{equation} \label{markovCVaRtransmission}
\inf_{\sigma_o^t} \mathds{E}_{\xi_t}\bigg\{[\mathcal{C}_o^t(D_t,e_t,\xi_t) + \sigma_o^t]_+  - \sigma_o^t \beta \bigg\} \leq 0,
\end{equation}
where $\sigma_o^t \in \mathbf{R}$ is an auxiliary variable \cite{rockafellar2000optimization}. The expression inside the expectation \eqref{markovCVaRtransmission} is expressed in the form

\begin{equation}\label{TransmissionCVaRConvex}
\underline{\mathcal{Q}}_o^t = \max_{k=1,2} \bigg[\langle \underline{\mathbf{a}}_{ok}(\mathbf{y}_t),\xi_t \rangle + \underline{\mathbf{b}}_{ok}(\sigma_o^t)\bigg],
\end{equation}
where the decision variable vector $\mathbf{y}$ consists of optimization variables $\{D, e, \sigma\}$. The convex expression \eqref{TransmissionCVaRConvex} is maximum of two affine functions, and we consider the following distributionally robust total CVaR objective
\begin{equation}\nonumber
\begin{split}
 \hat{J}_\textrm{Risk}^t = \sum_{\tau = t}^{t+\mathcal{H}_t}\sum_{o=1}^{N_\ell} \sup_{\mathds{Q}_\tau \in \hat{\mathcal{P}}_\tau^{N_s}}\mathds{E}^{\mathds{Q}_\tau} \max_{k=1,2} \bigg[ \langle \underline{\mathbf{a}}_{ok}(\mathbf{y}_\tau), \hat{\xi}_\tau \rangle + \underline{\mathbf{b}}_{ok}\sigma_o^\tau) \bigg].
\end{split}
\end{equation}
The above multi-period distributionally robust optimization can be equivalently reformulated as a quadratic program. The details of the linear reformulation are shown in \cite{mohajerin}.
The objective is to minimize a weighted sum of an operational cost function and the total worst-case CVaR of the affine constraints in the set $\mathcal{V}^t$ (e.g., line flow constraints and local device constraints).\\
\hspace{-2mm} \textbf{Data-based distributionally robust stochastic DC OPF}
\begin{subequations} \label{DROOPFTranmission}
\begin{align}
	&\inf_{\mathbf{y}_\tau, \sigma_o^\tau} \sum_{\tau = t}^{t + \mathcal{H}_t}\bigg\{\mathds{E}[\hat{J}_{\textrm{Cost}}^\tau] +  \rho \sup_{\mathds{Q}_\tau \in \hat{\mathcal{P}}_\tau^{N_\ell}}\sum_{o=1}^{N_\ell} \mathds{E}^{\mathds{Q}_\tau}[\underline{\mathcal{Q}_o^\tau}]  \bigg\}, \nonumber\\
	&=\inf_{\begin{subarray}{c}\mathbf{y}_\tau, \sigma_o^\tau,\\ \lambda_o^\tau, s_{io}^\tau, \varsigma_{iko}^\tau \end{subarray}} \sum_{\tau = t}^{t + \mathcal{H}_t} \bigg\{\mathds{E}[\hat{J}_\textrm{Cost}^\tau] + \sum_{o=1}^{N_\ell} \bigg(\lambda_o^\tau\varepsilon_\tau + \frac{1}{N_s}\sum_{i=1}^{N_s} s_{io}^\tau\bigg) \bigg\},\\
	&\textrm{subject to} ~~~~~~~~~~~~~~~~~~~~~~~~~~~~~~~~~~~~~~~~~~~~~~~~~~~~~~~~~~~~~ \nonumber\\
	&\sum_{d=1}^{N_d}\Big[ (r_d^{\bar{t}} + C_d^{\bar{t}}(A^d_{\bar{t}}x_{0}^{d}+B^d_{\bar{t}}e_{\bar{t}}^d))\Big]_{[\underline{t},\bar{t}]} = 0,\\
	& \sum_{d=1}^{N_d}\Big[(G_d^{\bar{t}} + C_d^{\bar{t}}B^d_{\bar{t}}D^d_{\bar{t}})\Big]_{[\underline{t},\bar{t}]}=0, \\
	& \rho (\underline{\mathbf{b}}_{ok}(\sigma_o^\tau) + \langle \underline{\mathbf{a}}_{ok}(\mathbf{y}_\tau), \hat{\xi}^{i}_\tau\rangle + \langle\varsigma_{iko}^\tau, \mathbf{d}-H\hat{\xi}^{i}_\tau\rangle) \le s_{io}^\tau,\\
&\|H^\intercal\varsigma_{iko}^\tau-\rho \underline{\mathbf{a}}_{ok}(\mathbf{y}_\tau)\|_\infty \le \lambda_o^\tau,\\
& \varsigma_{iko}^\tau \geq 0,\\
& \forall i\le N_s, \forall o \le N_\ell, k=1,2, \tau = t,\ldots, \mathcal{H}_t, \nonumber 
\end{align}
\end{subequations}
where $\rho \in \mathbf{R}_+$ quantifies the power system operators' risk aversion. Once again, this is a quadratic program that explicitly uses the training dataset $\hat{\Xi}^{N_s}_\tau := \{\hat{\xi}_\tau^{i}\}_{i\leq N_s}$, and the risk aversion parameter $\rho$ and the Wasserstein radius $\varepsilon_\tau$ allow us to explicitly balance tradeoffs between efficiency, risk, and sampling errors inherent in $\hat{\Xi}^{N_s}_\tau$. The uncertainty set is modeled as a polytope $\Xi := \{\xi \in \mathbf{R}^{N_\xi}: H\xi \leq \mathbf{d} \}$. The constraint $\varsigma_{iko} > 0$ holds since the uncertainty set is not-empty; on the other hand, in a case with no uncertainty (i.e, $\varsigma_{iko} = 0$), the variable $\lambda$  does not play any role and $s_{io}^\tau = \rho (\underline{\mathbf{b}}_{ok}(\sigma_o^\tau) + \langle \underline{\mathbf{a}}_{ok}(\mathbf{y}_\tau), \hat{\xi}_\tau^{i}\rangle)$.
\section{Conclusions}
In this paper, we have proposed a framework for the data-based distributionally robust stochastic OPF based on finite dataset descriptions of forecast error distributions across the power network. The method allows efficient computation of multi-stage feedback control policies that react to forecast errors and provide robustness to inherent sampling errors in the finite datasets. Tractability is obtained by exploiting convex reformulations of ambiguity sets based on Wasserstein balls centered on empirical distributions. The general framework is adapted to both distribution networks and transmission systems by allowing general device models and utilizing different linearizations of the AC power flow equations. The effectiveness and flexibility of our proposed method is demonstrated in a companion paper, Part II \cite{Guo2018DataDriven2}, which uses the method to handle a over-voltages in a distribution network with high solar penetration and to address line flow risks in a transmission system with high wind penetration. These numerical results also show that the method has superior out-of-sample performance and allows a more principled and systematic tradeoff of efficiency and risk.

\section*{Supplementary materials}
Implementation codes for 1) data-based distributionally robust stochastic OPF and 2) data-based distributionally robust stochastic DC OPF can be download from \href{https://github.com/TSummersLab/Distributionally-robust-stochastic-OPF}{[Link]}. Numerical case studies (e.g., overvoltage in distribution networks and $N$-1 security in transmission systems), and discussions on simulation results are presented in  \cite{Guo2018DataDriven2}.

\section*{Appendix}
Pertinent reformulations of the probabilistic line flow constraints in transmission systems are described in this section. Consider,  for simplicity, the time period between $[0,t]$ and rewrite the line flow constraint as
    \begin{equation}\nonumber
        \mathds{E}~\mathcal{R}\bigg\{\sum_{d=1}^{N_d}\Gamma_d^t(r_d^t + G_d^t\bm{\xi}_t + C_d^t\mathbf{x}_t^d) -\bar{p}_t \bigg\}\le 0,\\ 
    \end{equation}
    where $\mathbf{x}_t^d = A_t^dx_0^d + B_t^d\mathbf{u}_t^d$ and the affine reserve policies are $\mathbf{u}_t^d = D_t^d\bm{\xi}_t + e_t^d$. The CVaR counterpart of each individual constraints is given by [cf.~\eqref{markovCVaRtransmission}] 
    \begin{equation}\label{TransactionsLineAppendix}
    \begin{aligned}
     \inf_{\sigma_o^t} \mathds{E}_{\xi_t}\textrm{max}\bigg\{\Big[ \sum_{d=1}^{N_d}\Gamma_d^t(r_d^t + G_d^t\bm{\xi}_t + C_d^t(A_t^dx_0^d + B_t^dD_t^d\bm{\xi}_t &\\
     + B_t^de_t^d)) -\bar{p}_t\Big]_o + \sigma_o^t  \bigg\} - \sigma_o^t \beta  \leq 0 .&
     \end{aligned}
    \end{equation}
 Since the decision variables $\{D_t, e_t, \sigma_o^t\}$ enter as linear terms in \eqref{TransactionsLineAppendix}, a compact expression similar to \eqref{TransmissionCVaRConvex} can be obtained; i.e.,  
    \begin{equation}\nonumber
    \underline{\mathcal{Q}}_o^t = \max_{k=1,2} \bigg[\langle \underline{\mathbf{a}}_{ok}(\mathbf{y}_t),\xi_t \rangle + \underline{\mathbf{b}}_{ok}(\sigma_o^t)\bigg],
    \end{equation}
    where\\
~\\
$k = 1$,\\
\begin{equation*}
\underline{\mathbf{a}}_{o1}(\mathbf{y}_t)=\left[
\begin{matrix}
\big[\sum_{d=1}^{N_d}\Gamma_d^tC_d^tB_t^dD_t^d\big]_o,  \big[\sum_{d=1}^{N_d}\Gamma_d^tC_d^tB_d^te_t^d\big]_o
\end{matrix}\right],
\end{equation*}
\begin{equation*}
\underline{\mathbf{b}}_{o1}(\sigma_o^t) = \big[ -\bar{p} +\sum_{d=1}^{N_d}\Gamma_d^t(r_d^t + G_d^t\xi + C_d^tA_t^dx_0^d) \big]_o+ \sigma_o^t - \sigma_o^t\beta,\\
\end{equation*}
$k = 2$,\\
\begin{equation*}
\underline{\mathbf{a}}_{o2}(\mathbf{y}_t)=\left[
\begin{matrix}
\bm{0}_{N_\xi t}^\intercal, 0
\end{matrix}\right], 
\underline{\mathbf{b}}_{o2}(\sigma_o^t)= -\sigma_o^t\beta,
\end{equation*}
for all $o$. Similar steps can be followed to obtain the  constraints for each device  and the voltage magnitude.
\bibliography{Final_PartI_DROOPF}
\bibliographystyle{ieeetr}

\end{document}